\documentclass[10pt, a4paper, reqno]{amsart}

\usepackage{amsmath, amssymb, amsfonts, enumerate}
\usepackage[english]{babel}
\usepackage{fancyhdr}
\usepackage{graphics}
\usepackage{amsfonts}
\usepackage{mathrsfs}
\usepackage{mathptmx}
\usepackage{stmaryrd}
\usepackage{epsfig} 
\usepackage{hyperref}
\usepackage[utf8]{inputenc}

\theoremstyle{plain}
\newtheorem{Theorem}{Theorem}
\newtheorem{Proposition}{Proposition}
\newtheorem{Lemma}{Lemma}
\newtheorem{Corollary}{Corollary}
\newtheorem{Remark}{Remark}
\newtheorem{Example}{Example}
\newtheorem{Definition}{Definition}
\newtheorem{Problem}{Problem}

\makeatletter
\def\imod#1{\allowbreak\mkern10mu({\operator@font mod}\,\,#1)}
\makeatother

\begin{document}

\title[$\ $]{Left invertibility of I/O quantized
linear systems in dimension 1: a number theoretic approach}

\author[N. Dubbini]{Nevio Dubbini}
\address{Nevio Dubbini: Department of Mathematics ``L. Tonelli'',
University of Pisa, Pisa, Italy}
\email{dubbini@mail.dm.unipi.it}

\author[M. Monge]{Maurizio Monge}
\address{Scuola Normale Superiore, Pisa, Italy}

\author[A. Bicchi]{Antonio Bicchi}
\address{Interdepartmental center ``E. Piaggio'', University
of Pisa, Pisa, Italy}
\date{\today}

\begin{abstract}
This paper studies left invertibility of discrete-time linear I/O
quantized linear systems of dimension $1$. Quantized outputs are
generated according to a given partition of the state-space, while
inputs are sequences on a finite alphabet. Left invertibility, i.e.
injectivity of I/O map, is reduced to left D-invertibility, under
suitable conditions. While left invertibility takes into account
membership in sets of a given partition, left D-invertibility
considers only distances, and is very easy to detect. Considering
the system $x^+=ax+u$, our main result states that left
invertibility and left D-invertibility are equivalent, for all but a
(computable) set of $a$'s, discrete except for the possible presence
of two accumulation point. In other words, from a practical point of
view left invertibility and left D--invertibility are equivalent
except for a finite number of cases. The proof of this equivalence
involves some number theoretic techniques that have revealed a
mathematical problem important in itself. Finally, some examples are
presented to show the application of the proposed method.
\end{abstract}

\maketitle

\section{Introduction}

Left invertibility is an important problem of systems theory, which
corresponds to injectivity of I/O map. It deals with the possibility
of recovering unknown inputs applied to the system from the
knowledge of the outputs.

We investigate left invertibility of discrete--time linear I/O
quantized systems in a continuous state-space of dimension $1$. In
particular, inputs are arbitrary sequences of symbols in a finite
alphabet: each symbol is associated to an action on the system.
Information available on the system is represented by sequences of
output values, generated by the system evolution according to a
given partition of the state-space (uniform quantization).

In recent years there has been a considerable amount of work on
quantized control systems (see for instance \cite{Del}, \cite{Pica},
\cite{Szn} and references therein), stimulated also by the growing
number of applications involving ``networked'' control systems,
interconnected through channels of limited capacity (see e.g.
\cite{Bic,Car,Tat}). The quantization and the finite cardinality of
the input set occur in many communication and control systems.
Finite inputs arise because of the intrinsic nature of the actuator,
or in presence of a logical supervisor, while output quantization
may occur because of the digital nature of the sensor, or if data
need a digital transmission.

Applications of left invertibility include fault detection in
Supervisory Control and Data Acquisition (SCADA) systems, system
identification, and cryptography (\cite{Edel,Inou}). Invertibility
of linear systems is a well understood problem, first handled in
\cite{Brok}, and then considered with algebraic approaches (see e.g.
\cite{Silv}), frequency domain techniques (\cite{Mas}, \cite{Mas2}),
and geometric tools (cf. \cite{Mors}). Invertibility of nonlinear
systems is discussed in (\cite{Resp}). More recent work has
addressed the left invertibility for switched systems (\cite{Vu}),
and for I/O quantized contractive systems(\cite{Nev}).

The main intent of the paper is to show that the analysis of left
invertibility can be substituted, under suitable conditions, by an
analysis of a stronger notion, called left D-invertibility. While
left invertibility takes in account whether two states are in the
same element of a given partition, left D-invertibility considers
only the distance between the two states. For this reason left
D-invertibility is very easy to detect. For the system $x^+=ax+u$,
the condition under which left invertibility and left
D-invertibility are equivalent has to do with the existence of an
infinite (periodic) orbit inside a certain set and the contemporary
occurrence of an algebraic condition satisfied by $a$. This two
conditions are as a matter of fact not restrictive, and indeed the
main theorem (Theorem \ref{th_uldi_uli_finite}) states that the set
of $a$ such that left D-invertibility and left invertibility are not
equivalent is discrete but possibly $2$ accumulation points. In
other words from a practical point of view ULI and ULDI are
equivalent except for a finite number of cases (see Theorem
\ref{th_uldi_uli_finite}).$\ $

The main tools used in the paper are a generalization of a classical
density theorem of Kronecker, and some geometry of numbers. The
Kronecker's theorem has to do with density in the unit cube of the
fractional part of real numbers. By means of a particular
construction the problem of ``turning'' left D-invertibility into
left invertibility can be handled with a Kronecker-type density
theorem. Geometry of numbers helps us to show that, even if the
Kronecker's theorem has not a straightforward application (we do not
have density) we can obtain our result anyway (we have
$\epsilon-$density, with $\epsilon$ small enough).

The paper is organized as follows: section $2$ contains a precise
statement of the problem under study, while section $3$ concerns the
number theoretic background needed. Section $4$ shows the procedure
to prove the equivalence between left D-invertibility and left
invertibility: the rational case is treated first, to show in a more
direct way ideas involved. This section contains also the main
result of the paper (Theorem \ref{th_uldi_uli_finite}). In section
$5$ explicit calculations are done in a comprehensive example.
Conclusions and future work are explained in section $6$. Finally,
there is a ``special'' section, the $7$-th, in which we collect the
notations used in the paper.

\section{Statement of the problem}

\begin{Definition}\label{def_uniform_partition}
The uniform partition of rate $\delta$ of $\Bbb R$ is
\[\mathcal{P}=\left\{\mathcal{P}_{i}\right\}_{i\in\Bbb Z}=
 \left\{\ [i\delta,(i+1)\delta[ \right\}_{i\in\Bbb Z}.\ \diamondsuit\]
\end{Definition}

In this paper we consider discrete-time, time-invariant, I/O
quantized linear systems of the form
\begin{equation}\label{sistema_iniziale}
             \left\{\begin{array}{l}
             x(k+1)=ax(k)+bu(k)\\
             y(k)=q_{\mathcal{P}}\big(cx(k)\big)\\
             \end{array}\right.
\end{equation}
where $x(k)\in\Bbb{R}$ is the state, $y(k)\in\Bbb Z$ is the output,
$u(k)\in\mathcal{U}\subset\Bbb R$ is the input, and $a,b,c\in\Bbb
R$. The map $q_{\mathcal{P}}: \Bbb{R} \rightarrow \Bbb Z$ is induced
by the uniform partition $\mathcal{P} =
\left\{\mathcal{P}_i\right\}_{i\in\Bbb Z}$ of $\Bbb{R}$ of rate
$\delta$ through $q_{\mathcal{P}}:(x\in \mathcal P_i) \mapsto i$ and
will be referred to as the output quantizer. We assume that
$\mathcal{U}$ is a finite set of cardinality $n$.

\begin{Remark}
Without loss of generality in the system (\ref{sistema}) we can
suppose $\delta=1,b=1,c=1$.
\end{Remark}

\emph{Proof:} Operate the substitutions $x(k)=\frac{\delta}{c}x(k)$
and $u(k)=\frac{\delta}{cb}u(k)$. $\diamondsuit$\\
\\
So we consider only systems of the form
\begin{equation}\label{sistema}
    \left\{\begin{array}{l}
             x(k+1)=ax(k)+u(k)=f_{u(k)}(x(k))\\
             y(k)=\lfloor x(k)\rfloor,
           \end{array}
    \right.
\end{equation}
where $\lfloor\cdot\rfloor$ denotes the integer part. Indicate with
$f^{k_2}_{k_1}(x_0,u_1,\ldots,u_{k_2})$ the sequence of outputs
$(y_{k_1},\ldots,y_{k_2})$ generated by the system (\ref{sistema})
with initial condition $x_0$ and input string
$(u_1,\ldots,u_{k_2})$.

\begin{Definition}
A pair of input strings $\{u_i\}_{i\in\Bbb{N}}$,
$\{u'_i\}_{i\in\Bbb{N}}$ is uniformly distinguishable in $k$ steps,
$k\in\Bbb N$, (or with distinguishability time $k$) if there exists
$l\in\Bbb N$ such that $\forall (x_0, x'_0) \in \Bbb R^2$ and
$\forall m>l$ the following holds:
\[u_m\not=u'_m\;\Rightarrow\;f^{m+k}_m(x_0,u_1,\ldots,u_{m+k}) \not=
f^{m+k}_m(x'_0,u'_1,\ldots,u'_{m+k}).\] In this case, we say that
the strings are uniformly distinguishable with waiting time $l$.
$\diamondsuit$
\end{Definition}

\begin{Definition}
A system of type (\ref{sistema}) is uniformly left invertible (ULI)
in $k$ steps if every pair of distinct input sequences is uniformly
distinguishable in $k$ steps after a finite time $l$, where $k$ and
$l$ are constant. $\diamondsuit$
\end{Definition}

For a ULI system, it is possible to recover the input string until
instant $m$ observing the output string until instant $m+k$. For
applications, however it is important to obtain an algorithm to
reconstruct the input symbol used at time $m>l$ by processing the
output symbols from time $m$ to $m+k$.

\begin{Definition}\label{def_set_q}
Define
\[Q = \bigcup_{y\in\Bbb{Z}}\{q^{-1}(y)\times
q^{-1}(y)\} = \bigcup_{y\in\Bbb{Z}} \{[y,y+1[\times[y,y+1[\}
\subset\Bbb{R}^{2}\] i.e. the union of the preimages of two
identical output symbols. In other words, $Q$ contains all pairs of
states that are in the same element of the partition $\mathcal{P}$.
$\diamondsuit$
\end{Definition}

To address invertibility, we are interested in studying the
following system on $\Bbb{R}^{2}$:
\begin{equation}\label{sistema2d}
X(k+1) = F_{U(k)}(X(k)) = \left[\begin{array}{c}
f \left( x(k), u(k) \right) \\
f \left( x^{\prime}(k), u^{\prime}(k) \right)
\end{array}\right]
\end{equation}
where $X(k) = \left[\begin{array}{c} x(k) \\ x^{\prime}(k)
\end{array}\right], \ \ U(k) = \left(u(k), u^{\prime}(k) \right).$
If it is possible to find an initial state in $Q$ and an appropriate
choice of the strings $\{u_k\},\{u'_k\}$ such that the orbit of
(\ref{sistema2d}) remains in $Q$, it means that the two strings of
inputs give rise to the same output for the system (\ref{sistema}).
Therefore conditions ensuring that the state is outside $Q$ for some
$k$ will be seeked to guarantee left invertibility. We will need
another notion of left invertibility, stronger but very easy to
check, that we define in the following. It will be central in our
discussion.

\begin{Definition}\label{def_difference_system}
The difference system associated with the system (\ref{sistema}) is
\begin{equation}\label{equation_difference_set}
    z(k+1) = az(k)+v(k)=f_{v(k)}(z(k))
\end{equation}
where $z(k)\in\Bbb R$, $v(k)\in\mathcal{V}=\mathcal{U}-\mathcal{U}$.
$\diamondsuit$
\end{Definition}

\begin{Remark}\label{remark_difference_system_leq_1}
The difference system represents at any instant the difference
between the two states $z(k)=x(k)-x^{\prime}(k)$ when the input
symbols $u(k)-u^{\prime}(k)=v(k)$ are performed. So we are
interested in understanding the conditions under which
\[\{z(k)\}\ \cap\ ]-1,1[\ =\emptyset.\]
Indeed, this implies that $y(k)\not=y^{\prime}(k)$. The converse is
obviously not true. $\diamondsuit$
\end{Remark}

Indicate with $D^{k_2}_{k_1}(z_0,v_1,\ldots,v_{k_2})$ the sequence
$(z(k_1),\ldots,z(k_2))$ generated by the difference system with
initial condition $z_0$ and input string $(v_1,\ldots,v_{k_2})$.

\begin{Definition}
A pair of input strings $\{u_i\}_{i\in\Bbb{N}}$,
$\{u'_i\}_{i\in\Bbb{N}}$ is uniformly D-distinguishable in $k$
steps, $k\in\Bbb N$ (or with distinguishability time $k$), if there
exists $l\in\Bbb N$ such that $\forall (z_0) \in \Bbb R$ and
$\forall m>l$ the following holds:
\[v_m\not=0\ \Rightarrow\ D_m^{m+k}(z_0,v_1,\ldots,v_{m+k})\not\in\ ]-1,1[^{m+k+1},\]
where $v_i=u_i-u^{\prime}_i$. In this case, we say that the strings
are uniformly D-distinguishable with waiting time $l$.
$\diamondsuit$
\end{Definition}

\begin{Definition}
A system of type (\ref{sistema}) is uniformly left D-invertible
(ULDI) in $k$ steps if every pair of distinct input sequences is
uniformly D-distinguishable in $k$ steps after a finite time $l$,
where $k$ and $l$ are constant. $\diamondsuit$
\end{Definition}

\begin{Remark}
Thanks to Remark \ref{remark_difference_system_leq_1} uniform left
D-invertibility implies uniform left invertibility. $\diamondsuit$
\end{Remark}

\begin{Proposition}\label{prop_uldi_one_step}
The system (\ref{sistema}) is either ULDI in time 1, or not ULDI at
all, depending on the following condition is satisfied:
\[\min_{0\not=v\in\mathcal V} |v|\geq |a|+1.\]
\end{Proposition}

\emph{Proof:} A sufficient condition for uniform left
D-invertibility in one step is
\[\forall v\in\mathcal{V},v\not=0:\ |v|\geq |a|+1:\]
indeed in this hypothesis $\forall v\in\mathcal{V},v\not=0$
\[ ]-1,1[\ \cap\ \big\{a\cdot(]-1,1[)+v\big\} =\ ]-1,1[\ \cap\ ]-a+v,a+v[\ = \emptyset\]

We now prove that if $\exists v\in\mathcal{V},v\not=0:\ |v|< |a|+1$,
then the system is not uniformly left D-invertible. Indeed in this
case the system
\[\left\{\begin{array}{c}
    ax_1+v=x_2 \\
    ax_2-v=x_1
  \end{array}\right.\]
has the solution $x_1=\frac{-v}{a+1},x_2=\frac{v}{a+1}$. Since
$|x_1|,|x_2|<1$ the difference system has the infinite orbit
$\{x_1,x_2,x_1,x_2,\ldots\}\subset]-1,1[$. Therefore system
(\ref{sistema}) is not left D-invertible. $\diamondsuit$
$\quad$\\

Proposition \ref{prop_uldi_one_step} shows a trivial way to check
ULDI for systems (\ref{sistema}). The problem under study is the
following:

\begin{Problem}
State mathematical conditions for the equivalence between ULDI and
ULI of a uniformly quantized linear system of the form
(\ref{sistema_iniziale}). $\diamondsuit$
\end{Problem}

\section{Mathematical background}

We will mainly need results from number theory: our proofs are
essentially based on the application of a density Theorem of
Kronecker (see \cite{Har}), sufficient in the case in which $a$ is
trascendental. For the algebraic case we need further computations
involving the Mahler measure of polynomials.

\begin{Definition}
The numbers $\vartheta_1,\ldots,\vartheta_M\in\Bbb R$ are linearly
independent over $\Bbb Z$ if the following holds:
\[ k_1,\ldots,k_M\in\Bbb Z: k_1\vartheta_1+\ldots,+k_M\vartheta_M=0\ \ \Rightarrow\ \ k_1=\ldots=k_M=0.\ \diamondsuit\]
\end{Definition}

\begin{Theorem}[Kronecker]\cite{Har}\label{theorem_Kro_2}
If $\alpha_1,\ldots,\alpha_M,1\in\Bbb R$ are linearly independent
over $\Bbb Z$, then, for every
$\vartheta_1,\ldots,\vartheta_M\in\Bbb R$ the set of points
\[\left\{ \left[frac(l\alpha_1+\vartheta_1),\ldots,frac(l\alpha_M+\vartheta_M)\right]: l\in\Bbb R \right\}\]
is dense in the unit cube of $\Bbb R^M$. $\diamondsuit$
\end{Theorem}

\begin{Definition}
A set of independent linear relations among
$\alpha_1,\ldots,\alpha_M\in\Bbb R$ is said to be maximal if no
other independent linear relation can be found among these numbers.
$\diamondsuit$
\end{Definition}

\begin{Remark}\label{remark_kron}
A corollary of the Kronecker's Theorem (clear from the proof) is
that, if the numbers $\alpha_1,\ldots,\alpha_M\in\Bbb R$ satisfy a
maximal set of nontrivial linear equations $L =
\{L^j(\alpha_1,\ldots,\alpha_M)=0\ \ for\ j=1,\ldots,J\},$ then the
set of points
\[\left\{
\left[frac(\alpha_1+l\vartheta_1),\ldots,frac(\alpha_M+l\vartheta_M)\right]:
l\in\Bbb R \right\}\] is dense in
\[C_L = frac(\{x_1,\ldots,x_M:\ L^j(x_1,\ldots,x_M)=0,\ j=1\ldots,J\}).\ \diamondsuit\]
\end{Remark}

\begin{Definition}
We define the set of linear relations $L$ to be integer-maximal for
the numbers $\alpha_1,\ldots,\alpha_M$ if
\begin{itemize}
  \item The set of points $\left\{\left[frac(\alpha_1+l\vartheta_1),\ldots,frac(\alpha_M+l\vartheta_M)\right]:
l\in\Bbb R \right\}$ is dense in $C_L$;
  \item The linear relations $L^j$ are formed with integer
  coefficients;
  \item $C_L = \{x\in\Bbb R^M:\ L^jx\in\Bbb Z^J\}\ \cap\ [0,1[^M$.
  $\diamondsuit$
\end{itemize}
\end{Definition}

\begin{Definition}
A number $\rho\in\Bbb C$ is called algebraic if there exists a
polynomial $R(x)\in\Bbb Z[x]$ such that $R(\rho)=0$.  In this case
there exists a unique monic polynomial $R(x)\in\Bbb Z[x]$ with
minimal degree $q$. $R(x)$ is called the minimal polynomial of
$\rho$ and $q$ its degree. A number $\rho\in\Bbb C$ is called
trascendental if it is not algebraic. $\diamondsuit$
\end{Definition}
Note that $1,\rho,\rho^2,\ldots,\rho^M$ are linearly independent if
and only if the degree of $\rho$ is at least $M+1$.

\begin{Definition}
The $i-th$ symmetric polynomial in $q$ variables is
\[e_i(x_1,\ldots,x_q) = \sum_{1\leq j_1\leq\ldots\leq j_q} x_{j_1}\cdot\ldots\cdot x_{j_q}.\]
\end{Definition}

\begin{Definition}\label{def_mahler_measure}
If $R(x)$ is the polynomial
\[ R(x) = \sum_{i=0}^q r_i x^i = r_q \cdot \prod_{j=1}^q (x-\rho_j),\]
where the $\rho_j$'s are the roots of the polynomial, its Mahler
measure is defined as
\[ \mathfrak M(R) = r_q \cdot \prod_{j=1}^q \max\left\{1,|\rho_j|\right\}.\]
\end{Definition}

Mahler measure has many interesting properties. For instance, since
$r_i$ is equal to $r_q$ multiplied by the $i$-th symmetric
polynomial of the $\rho_j$, which is made of precisely
$\binom{q}{i}$ monomials in the $\rho_i$ where each $\rho_i$ appears
with degree at most $1$, we have that $r_i$ is sum of $\binom{q}{i}$
terms each $\leq \mathfrak M(A)$ in absolute value, and consequently
\begin{equation} \label{mahlest}
\left|r_i\right| \leq \binom{q}{i} \cdot \mathfrak
M(R),\qquad\text{for }0\leq i \leq q.
\end{equation}

If $\|R\|_\infty$ is the norm
\begin{equation}\label{eq_norm_infty}
\|R\|_\infty = \max_{0\leq i \leq q} \left|r_i\right|,
\end{equation}
we obtain from \eqref{mahlest}
\begin{equation}\label{mahlest2}
  \|R\|_\infty \leq \binom{q}{[q/2]} \cdot \mathfrak M(R).
\end{equation}$\ $\\
In the following we will also have to consider the quantity
\[ \mathfrak M(R(x/2)) = r_q \cdot \prod_{j=1}^q \max\left\{\frac{1}{2},|\rho_j|\right\}, \]
i.e. the Mahler measure of the polynomial $R(x/2)$. The last
equality is easily proved since
\begin{equation}
  \mathfrak M(R(x/2)) = \frac{r_q}{2^q} \cdot \prod_{j=1}^q \max\left\{1,|2\rho_j|\right\}\\
  = r_q \cdot \prod_{j=1}^q \max\left\{\frac{1}{2},|\rho_j|\right\}.
\end{equation}
In particular, note that
\begin{equation} \label{mahlest3}
  \mathfrak M(R) \leq 2^q \cdot \mathfrak M(R(x/2)).
\end{equation}

\section{ULI: the number theoretic approach}

Our strategy is the following: for trascendental $a$ in the system
(\ref{sistema}), we prove that ULDI is equivalent to ULI. Moreover,
for $a$ algebraic (and rational), we will show that these two
notions are very close, in a sense precisely specified
later.\\

\textbf{Notations:} Consider the system of dimension $2$ given by
(\ref{sistema2d}), and suppose that there exists at least one proper
orbit included in the set
\begin{equation}\label{eq_set_qprime}
Q^{\prime} = \Big\{\left(
            \begin{array}{c}
                t \\
                t \\
            \end{array}
            \right) + \left(
                        \begin{array}{c}
                            s \\
                            0 \\
                        \end{array}
                    \right):\ s\in\ ]-1,1[\ ,\ t\in\Bbb R\Big\}.
\end{equation}
(such an orbit exists if and only if system (\ref{sistema}) is
\emph{not} ULDI). Take as initial condition $X_t(0) = \left(
                                                   \begin{array}{c}
                                                     t \\
                                                     t \\
                                                   \end{array}
                                                 \right)+\left(
                                                           \begin{array}{c}
                                                             s \\
                                                             0 \\
                                                           \end{array}
                                                         \right)\in\Bbb R^2$, with $t$, considered as a parameter,
varying in $\Bbb R$ and $s\ \in]-1,1[$ fixed. Then, for fixed input
string
\begin{equation}\label{eq_traj_in_q_prime}
X_t(k)=\left(
          \begin{array}{c}
            a^kt+a^ks+a^{k-1}u_1+\ldots+u_k \\
            a^kt+a^{k-1}u^{\prime}_1+\ldots+u^{\prime}_k \\
          \end{array}
        \right) = a^k\left(
                    \begin{array}{c}
                      t \\
                      t\\
                    \end{array}
                  \right)+\left(
                            \begin{array}{c}
                              c_k \\
                              c^{\prime}_k \\
                            \end{array}
                          \right).\end{equation}
Suppose that an orbit $\{X_t(i)\}_{i=1}^{\infty}$ is included in
$Q^{\prime}$. We can see the points $X_t(i)$, when $t$ varies in
$\Bbb R$, as points moving along the line
\begin{equation}\label{eq_varrho}
\varrho_i = \Big\{a^i\left(\begin{array}{c}
                      t \\
                      t\\
                    \end{array}
                  \right)+\left(
                            \begin{array}{c}
                              c_i \\
                              c^{\prime}_i \\
                            \end{array}\right):\ t\in\Bbb R\Big\}
\end{equation}
with initial condition $\left(\begin{array}{c}
                              c_i \\
                              c^{\prime}_i \\
                            \end{array}\right)$ and velocity $a^i$.
Call $k_i$ the distance between the point $\left(\begin{array}{c}
                                                        c_i \\
                                                        c^{\prime}_i \\
                                                    \end{array}\right)$
and the union of positive coordinate axes along the line $\varrho_i$
(refer to the figure $1$).
\begin{figure}
\begin{center}
\includegraphics[width=0.7\textwidth]{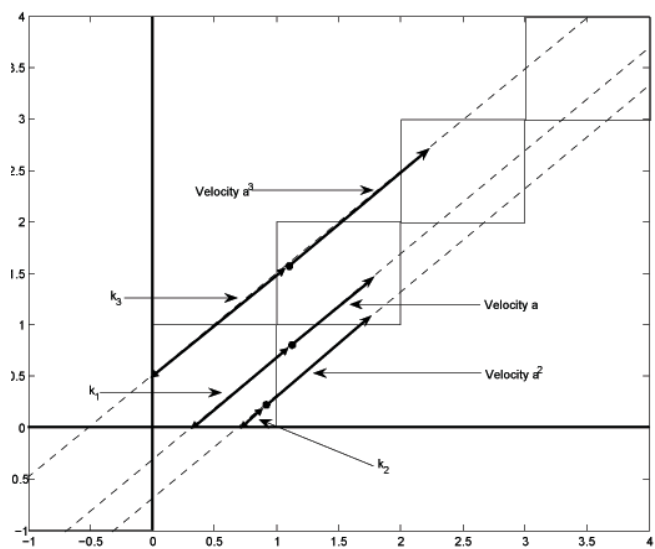}
\caption{Here, for $i=1,2,3$, the point $X_t(i)$ has a distance
$k_i$ from the union of positive coordinate axes along the line
$r_i$ (drawn with a dashed line), and ``velocity'' $a^i$ (with
respect to $t$)}
\end{center}
\end{figure}\\

\subsection{Trascendental $a$}
$\ $\\

The following technical lemma gives a necessary condition for
uniform left invertibility, a basilar ingredient in the proof of
Theorem \ref{theorem_a_trascendental}.

\begin{Lemma}\label{lemma_ij_kron_cond}
Consider the $2$-dimensional system (\ref{sistema2d}) with the
notations just introduced. Suppose that $\forall\epsilon>0$,
$\forall J\in\Bbb N$, $\forall s\in]-1,1[$, $\forall
\{U(j)\}_{j\in\Bbb N}$ there exists $t\in\Bbb R$ such that, if
$\{X(j)\}_{j=0}^J\subset Q^{\prime}$ is the orbit give by
$X(0)=\left(
                                                                  \begin{array}{c}
                                                                    t+s \\
                                                                    t \\
                                                                  \end{array}
                                                                \right)$ and input sequence
$U(j)$, the following holds for every $j=1,\ldots,J$:
\begin{equation}\label{eq_ij_kron_cond}
frac\left( k_j+ta^j \right)<\epsilon.
\end{equation}
Then the system is not ULI.
\end{Lemma}

\emph{Proof:} Suppose that an orbit $\{X_{t}(j)\}_{j=1}^{\infty}$ is
included in $Q^\prime$. Observe that $frac\left( k_j+a^jt \right)=0$
if and only if $X_t(j)$ belongs to some translation of
\begin{equation}\label{eq_set_omega}
\Omega=[0,1]\times\{0\}\ \cup\ \{0\}\times[0,1]
\end{equation}
along the diagonal of $\Bbb R^2$, that is entirely included in $Q$,
i.e. a translation that takes $\Omega$ to the ``bottom-left
boundary'' of a square of $Q$. It's now easy to see that, for every
$X_t(j)$ there exists $\epsilon>0$ such that, if $frac\left(
k_j+a^jt \right)<\epsilon$ then $X_t(j)\in Q$. Therefore, if the
relations (\ref{eq_ij_kron_cond}) are satisfied, then there exists
an arbitrary long orbit included in $Q$. $\diamondsuit$

\begin{Proposition}\label{prop_unid_algebraic_a}
Suppose that the system (\ref{sistema}) is not ULDI. If $a$ is an
algebraic number of degree $K$ then the system is not ULI in $K-1$
steps.
\end{Proposition}

\emph{Proof:} Since the system is not ULDI there exist arbitrary
long orbits included in $Q^{\prime}$. Fix one of these orbits of
length greater than $K-1$.

If, for every $\epsilon>0$, and every $k_1,\ldots,k_{K}\in\Bbb R$
there exists a $t\in\Bbb R$ such that
\begin{equation}\label{eq_kron_cond}
frac \left(k_i+a^it \right)<\epsilon\ for\ i=0,\ldots,K-1
\end{equation}
then the system (\ref{sistema}) is not ULI in $K-1$ steps by Lemma
\ref{lemma_ij_kron_cond}. Equation (\ref{eq_kron_cond}) is
equivalent to find integers $N_0,\ldots,N_K$ such that for every
$i=0,\ldots,K$
\[N_i\leq k_i+a^it<N_i+\epsilon\]
But, if $a$ is algebraic of degree $K$, then numbers $a^i$,
$i=0,\ldots,K-1$ are linearly independent over $\Bbb Z$, and by
Theorem \ref{theorem_kro_2} there always exists a $t$ such that
equation (\ref{eq_kron_cond}) holds, and so the system is not
uniformly left invertible in $K-1$ steps. $\diamondsuit$\\
\\
The following Theorem can be deduced immediately from Proposition
\ref{prop_unid_algebraic_a}.

\begin{Theorem}\label{theorem_a_trascendental}
Suppose that $a$ is trascendental. Then the system (\ref{sistema})
is ULI if and only if it is ULDI.
\end{Theorem}

\emph{Proof:} Suppose that system (\ref{sistema}) is not ULDI,
Proposition \ref{prop_unid_algebraic_a} states that, if the system
is ULI in $K$ steps then $a$ cannot be algebraic of degree greater
than $K+1$. The result follows easily since a trascendental number
is not algebraic of any degree. $\diamondsuit$

\begin{Corollary}
Consider the unidimensional system (\ref{sistema}), with
trascendental $a$. Then it is either ULI in one step, or it is not
ULI. $\diamondsuit$
\end{Corollary}$\ $\\

\subsection{Algebraic $a$}
$\ $\\

Suppose now that $a$ is algebraic of degree $K$, and that the
minimum polynomial of $a$ is
\[\alpha_Kt^{K-1}+\ldots+\alpha_0\in\Bbb Z[t].\]
We are interested in finding an $\epsilon$ (the minimum $\epsilon$)
such that for every $J\in\Bbb N$ there exists an $i\in\Bbb N$ and a
point in
\[[0,\epsilon[^{i+J+1}\]
in the sequence
\[\left\{frac\left(k_i+a^it\right),\ldots,frac\left(k_{i+J}+a^{i+J}t\right)\right\}\]
for every orbit $\{X_k\}_{k\in\Bbb N}\subset Q^{\prime}$.
Considering the $J-$dimensional torus $T^J$, the linear manifold (of
dimension $K$) associated with $a$, i.e. the linear manifold whose
image mod $1$ is what we called $C_L$ in Remark \ref{remark_kron},
is given by the following equations:
\begin{equation}\label{eq_pkj}
P_K^{J}\ =\ (t_0,\ldots,t_{J+K})\in\Bbb R^{J+K+1}:\
\left\{\begin{array}{l}
  \alpha_0t_0+\alpha_1t_1+\ldots+\alpha_Kt_K=0 \\
  \alpha_0t_1+\alpha_1t_2+\ldots+\alpha_Kt_{K+1}=0 \\
  \vdots \\
  \alpha_0t_{J}+\alpha_1t_{J+1}+\ldots+\alpha_Kt_{J+K}=0 \\
\end{array}\right.
\end{equation}

\begin{equation}\label{eq_psikj}
=\ Ker(\Psi_K^J) = Ker \left\{\left(
                        \begin{array}{cccccccc}
                            \alpha_0 & \alpha_1 & \ldots & \alpha_K & 0 & \ldots & 0 \\
                            0 & \alpha_0 & \alpha_1 & \ldots & \alpha_K & \ddots & 0 \\
                            \vdots &  \ddots & \ddots &  &  & \ddots & 0\\
                            0 & \ldots & 0 & \alpha_0 & \ldots & \alpha_{K-1} & \alpha_K\\
                        \end{array}
                        \right)\right\}.
\end{equation}
In other words, $frac(P_K^J)$ is the set in which the sequence
$\left\{frac\left(k_i+a^it\right),\ldots,frac\left(k_{i+J}+a^{i+J}t\right)\right\}$
is dense (by Remark \ref{remark_kron}).

\begin{Definition}\label{def_eps_a}
Denote with $\mathbf{\beta}_J$ the vector $(a,\ldots,a^J)$, and
define
\begin{equation}\label{eq_eps_a_1}
\epsilon(a) = \sup_{\mathbf{\zeta}\in\Bbb R^{J}}\ \ \sup_{J\in\Bbb
N}\ \ \inf_{t\in\Bbb R}\ \
\max_{l\in[i,i+J]}[frac(\mathbf{\zeta}+t\mathbf{\beta}_J)]_l\ =
\end{equation}
\begin{equation}\label{eq_eps_a_2}
=\ \inf \left\{ \epsilon \in\Bbb R \ :\  [0,\epsilon]^{K+J+1} \cap
\left(P_K^J+v+\Bbb Z^{K+J+1}\right) \neq \emptyset,\quad\text{for
each }v\in\Bbb R^{K+J+1},J\in\Bbb N\right\}
\end{equation}
\end{Definition}

Let us explain the meaning of $\epsilon(a)$. Suppose we are given
\emph{any} trajectory of the 2-dimensional system (\ref{sistema2d})
included in $Q^\prime$. Then, letting $t$ vary as a parameter, it
has the form (\ref{eq_traj_in_q_prime}), and we can investigate ULI
looking at fractional parts of $k_i+ta^i$, for $i=1,\ldots,J$, for
every $J\in\Bbb N$. Now, \emph{modulo} the $k_i$'s (i.e. modulo the
inputs), that is taking the $\sup$ on $\zeta$ in the definition,
$\epsilon(a)$ is the smallest $\epsilon$ such that for every
$J\in\Bbb N$ there exists $t\in\Bbb R$
\[frac(k_i+ta^i)<\epsilon\ \ for\ i=1\ldots,J.\]
It's now easy to see, looking at Lemma \ref{lemma_ij_kron_cond} that
$\epsilon(a)$ is useful to put in relation ULDI with ULI. Moreover,
by Remark \ref{remark_kron}, the set
\[\left\{frac(ta^i):\ t\in\Bbb R,\ i=1,\ldots,J\right\}\]
is dense in $frac(P_K^J)$. So $\epsilon(a)$ equals the following
quantity:
\[\epsilon(a) = \sup_{\mathbf{\zeta}\in\Bbb R^{J}}\ \ \sup_{J\in\Bbb N}\ \max_{l\in[i,i+J]} \left[frac(\mathbf{\zeta}+P_K^J)\right]_l.\]
This implies the second equivalent definition of $\epsilon(a)$ in
the definition \ref{def_eps_a}.

\begin{Proposition}
The map defined by the matrix $\Psi_K^J : \Bbb Z^{J+K+1} \rightarrow
\Bbb Z^{J+1}$ is surjective.
\end{Proposition}

\emph{Proof:} This is an immediate consequence of \cite[Lemma 2,
Chap. 1]{Cas}, which says that a rectangular integer $m\times l$
matrix, for $m>l$, can be completed to a square invertible integer
$m\times m$ matrix with determinant $1$ if and only if the greatest
common divisor of the $l\times l$ minors is $1$. Now, if a
rectangular integer $m\times l$ matrix, for $m>l$, can be completed
to an invertible integer $m\times m$ matrix, then the original
matrix must be clearly surjective from $\Bbb Z^m \rightarrow \Bbb
Z^l$.

All we have to do to apply the lemma is checking that the greatest
common divisor of the $k\times k$ minors of $\Psi_K^J$ is $1$, but
this is easy since for each prime $p$ we can consider the first
coefficient $\alpha_i$ of our polynomial such that $p$ does not
divide $\alpha_i$, and take the $k\times k$ minor made of the
columns $K-i+1, K-i+2, \dots, K-i+k$. Since this minor is lower
triangular when reduced modulo $p$ with all the elements on the
diagonal equal to $\alpha_i \imod p$, its determinant does not
vanish modulo $p$, and we are done. $\diamondsuit$

\begin{Proposition}\label{prop_eps_a}
\[\epsilon(a) = \inf \left\{ \epsilon \in \Bbb R \ :\
\left(\Psi_K^J\cdot[0,\epsilon]^{K+J+1}+w\right)
  \cap \Bbb Z^{J+1} \neq \emptyset,\quad\text{for each }w\in\Bbb
  R^{J+1}\right\}\]
\end{Proposition}

\emph{Proof:} First note that, for the set of vectors such that
$\Psi_K^J \cdot v$ has integer components, it holds
\[S = \left\{ v \in \Bbb R^{K+J+1}\ :\ \Psi_K^J \cdot v \in \Bbb Z^{J+1} \right\}\ =\ Ker\left(\Psi_K^J\right) + \Bbb Z^{K+J+1}.\]
Indeed, $Ker \Psi_K^J + \Bbb Z^{K+J+1} \subseteq S$ clearly, and for
each vector $w \in S$ there exist a vector $z \in \Bbb Z^{K+J+1}$
such that $\Psi_K^J \cdot w = \Psi_K^J \cdot z$, and consequently
the difference $v = w-z$ is in $Ker \left(\Psi_K^J\right)$, and we
have that $w = v+z \in Ker \Psi_K^J+\Bbb Z^{K+J+1}$. Now, since in
the (\ref{eq_eps_a_2}) we are quantifying over all vectors $v \in
\Bbb R^{K+J+1}$, we can equivalently say that
\begin{align} \label{prob2}
  \epsilon(a) = \inf \left\{ \epsilon \in \Bbb R \ :\  \left([0,\epsilon]^{K+J+1}+v\right)
  \cap S \neq \emptyset,\quad\text{for each }v\in\Bbb R^{K+J+1}\right\} \notag\\
\end{align}
applying the matrix $\Psi_K^J$ to the expression, and where we
denoted $\Psi_K^J\cdot[0,\epsilon]^{K+J+1}$ the image of
$[0,\epsilon]^{K+J+1}$ under the map $\Psi_K^J$. This passage must
be justified because the matrix $\Psi_K^J$ clearly does not have
rank $K+J+1$, but since $S$ contains \emph{all} the vectors that are
mapped to $\Bbb Z^{J+1}$ the first intersection will be non-empty
whenever the second one is (the other direction being trivial).
$\diamondsuit$\\

Following Proposition \ref{prop_eps_a}, we are investigating how big
must be $\epsilon$ to ensure that each set obtained translating
$\Psi_K^J\cdot[0,\epsilon]^{K+J+1}$ contains an integer vector. This
will be true if and only if
\[\Psi_K^J\cdot[0,\epsilon]^{K+J+1} + \Bbb Z^{J+1}=\Bbb R^J+1,\]
and equivalently if and only if $\Psi_K^J\cdot[0,\epsilon]^{K+J+1}$
contains a representative for each class in $\Bbb R^{J+1}/\Bbb
Z^{J+1}$.

\begin{Theorem}\label{th_notuldi_notuli}
Suppose that in the system (\ref{sistema}) there exists an infinite
orbit of the difference system in $]-1+\epsilon(a),1-\epsilon(a)[$.
Then the system is not uniformly left invertible.
\end{Theorem}

\emph{Proof:} In the hypotheses of the Theorem we can find, for
every $J\in\Bbb N$, an orbit of the 2-dimensional system
(\ref{sistema2d}) $X(0),\ldots,X(J)$, such that for every $i\in
0,\ldots,J$
\[
X(i)\in\left\{\left(
                  \begin{array}{c}
                    t \\
                    t \\
                  \end{array}
                \right)+ \left(
                           \begin{array}{c}
                             s \\
                             0 \\
                           \end{array}
                         \right):\ t\in\Bbb R,\ s\in [-1+\epsilon(a),1-\epsilon(a)]\right\},
\]
\[
frac(k_i)<\epsilon(a).
\]
This implies clearly that the system is not ULI (see figure $2$).
$\diamondsuit$

\subsubsection{$a\in\Bbb Q$}
$\quad$\\

We investigate first the rational case, because the estimates are
easier, and the results are straightforward. For the algebraic case
we need an harder work. Suppose that $a=\frac{p}{q}\in\Bbb Q$, with
$gcd(p,q)=1$. Then the minimal polynomial of $a$ is $P_a(x) = qx-p$.
So:

\begin{equation}\label{eq_associated_with_a_in_tk}
P_1^J = \left\{\begin{array}{c}
                pt_0+qt_1=0 \\
                pt_1+qt_2=0 \\
                \vdots \\
                pt_{J-2}+qt_{J-1}=0 \\
                pt_{J-1}+qt_J=0
                \end{array}\right.
\end{equation}

\begin{Proposition}\label{prop_estimate_rational}
Suppose that in the system (\ref{sistema}) $a=\frac{p}{q}\in\Bbb Q$.
Then $\epsilon(a)\leq \min\{\frac{1}{p},\frac{1}{q}\}$.
\end{Proposition}

\emph{Proof:} We show that the image of cube
$\left[0,\frac{1}{p}\right]^{J+2}$ under $P_1^J$ assumes each value
modulo $\Bbb Z^{J+1}$, and this can easily be done inductively in
the following way. Let $w = (w_1,\ldots,w_{J+2}) \in \Bbb R^{J+2}$,
we will build a vector $v = (v_1,\ldots,v_{J+1})$ such that
$\Psi_k^J\cdot w - v \in \Bbb Z^{J+1}$. Suppose that $v$ is such
that the first $i>0$ components of $\Psi_K^J\cdot w - v$ are in
$\Bbb Z$, and observe that while $w_{i+1}$ varies in the interval
$[0,1/p]$ the $i+1$-th component of $P_K^J\cdot w$ varies in an
interval large $1$, while the first $i$ components of $P_K^J\cdot w$
stay fixed. Consequently we can change $w_{i+1}$ to ensure that the
first $i+1$ components of $P_K^J\cdot w - v$ are in $\Bbb Z$, and
continuing in this way we prove our assertion. If $q \geq p$ we can
clearly proceed similarly but downwards, starting from the last
component. $\diamondsuit$

\begin{Corollary}
Suppose that in the system (\ref{sistema}) there exists an infinite
orbit of the difference system in
\[ \left]-1+\min\left\{\frac{1}{p},\frac{1}{q}\right\},1-\min\left\{\frac{1}{p},\frac{1}{q}\right\}\right[\]
Then the system is not uniformly left invertible. $\diamondsuit$
\end{Corollary}
\begin{figure}
\begin{center}
\includegraphics[width=0.9\textwidth]{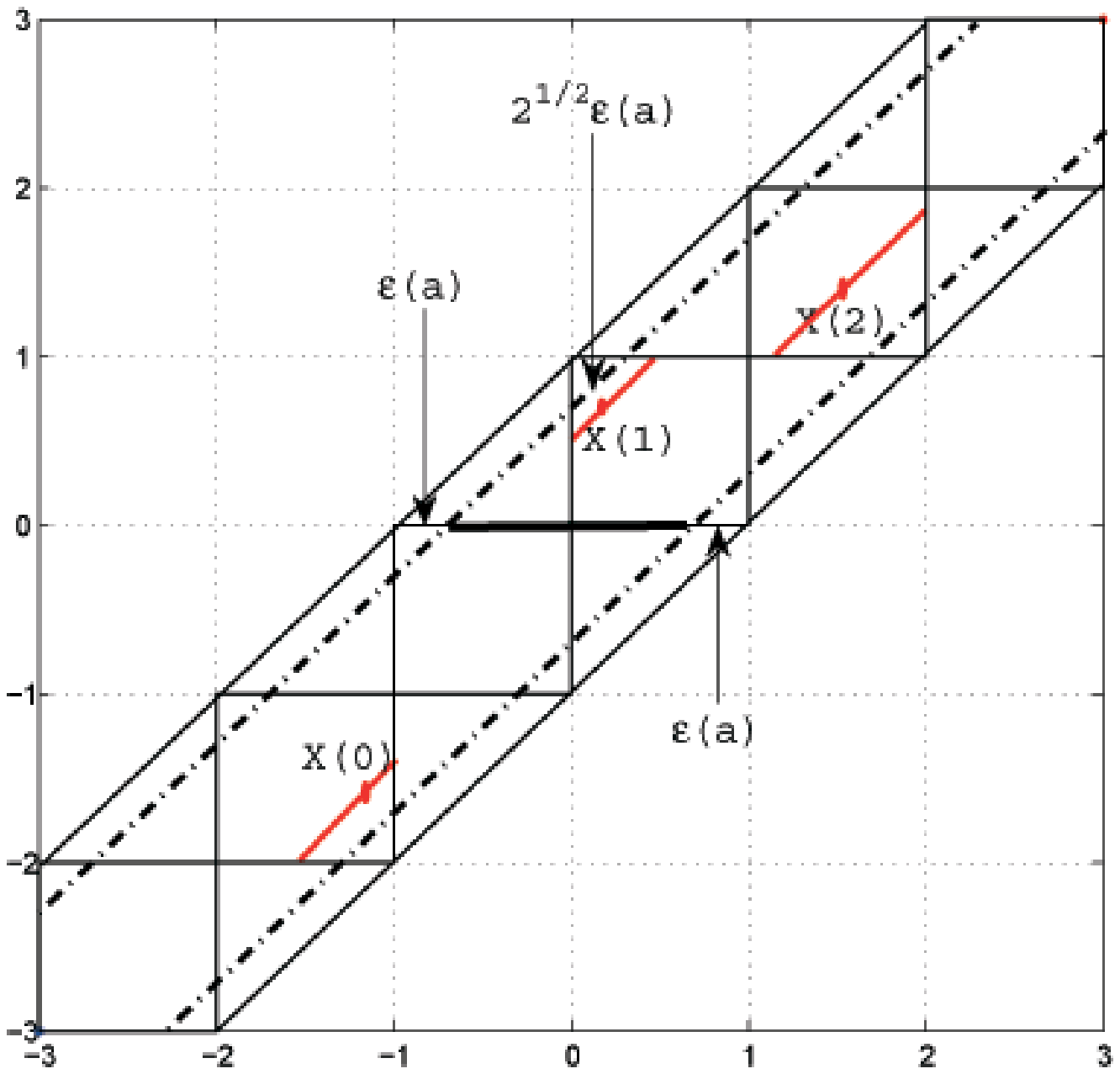}
\caption{In the hypotheses of Theorem \ref{th_notuldi_notuli} we can
find a trajectory (we represent here only $X(0),X(1),X(2)$) inside
the ``strip'' $\{[-1+\epsilon(a),1-\epsilon(a)]+ (t,t):\ t\in\Bbb
R\}$, drawn inside the dashed-dotted line.}
\end{center}
\end{figure}

\subsection{Algebraic $a\not\in\Bbb Q$}
$\quad$\\

Suppose that $a$ is algebraic, with minimum polynomial
$P_a(x)=\alpha_0+\ldots+\alpha_Kx^K$ of degree $K>2$. Then denote
with

\[
P_K^{J}\ =\left\{ (t_0,\ldots,t_{J+K})\in\Bbb R^{J+K+1}:\
\left\{\begin{array}{l}
  \alpha_0t_0+\alpha_1t_1+\ldots+\alpha_Kt_K=0 \\
  \alpha_0t_1+\alpha_1t_2+\ldots+\alpha_Kt_{K+1}=0 \\
  \vdots \\
  \alpha_0t_{J}+\alpha_1t_{J+1}+\ldots+\alpha_Kt_{J+K}=0 \\
\end{array}\right.\right\} = Ker\left(\Psi_K^J\right).
\]

\begin{Proposition}\label{prop_mahler}
Indicating with $P_a(x)$ the minimal polynomial of an algebraic
number $a$, the following estimate holds:
\begin{equation}\label{eq_estimate_mahler}
\epsilon(a) \leq \min\left\{ \frac{1}{\mathfrak M(P_a(x/2))},
\frac{1}{\mathfrak M(2^{-K}P_a(2x))}\right\}.
\end{equation}
\end{Proposition}

\emph{Proof:} See Theorem $1$ of \cite{Dub2}. $\diamondsuit$

\begin{Theorem}\label{theo_epsilon_a_alg}
Indicating with $P_a(x)=\alpha_Kt^k+\ldots+\alpha_0$ the minimal
polynomial of an algebraic number $a$, the following estimate holds:
\[\epsilon(a)\leq const\cdot\min\left\{\frac{1}{|\alpha_i|}:\
i=1,\ldots,k\right\},\] with the constant depending only on the
degree of $a$. Moreover, the constant is less than or equal to
\[\min_i \left\{ \binom{k}{i}2^{\min\ \{i,k-i\}}
\frac{1}{|\alpha_i|} \right\}.\]
\end{Theorem}

\emph{Proof:} By Proposition \ref{prop_mahler} it holds the estimate
(\ref{eq_estimate_mahler}). The two terms $2^{-K}\mathfrak
M(P_a(2x))$ e $\mathfrak M(P_a(x/2))$ are the Mahler measures of
respectively the polynomial with coefficients
$2^{-K}\alpha_k,\dots,2^{-1}\alpha_1,\alpha_0$, and the polynomial
with coefficients $\alpha_k,\dots,2^{-k+1}\alpha_1,2^{-k}\alpha_0$.
Moreover by (\ref{mahlest2}) it holds
\[2^{-i}\alpha_i \leq \binom{k}{i} M(P(x/2)),\]
\[2^{-k+i}\alpha_i \leq \binom{k}{i} M(2^{-k}P(2x))\]
and we are done. $\diamondsuit$

\begin{Corollary}
Suppose that in the system (\ref{sistema}) the degree of $a$ is at
least $2$. Suppose that there exists a proper path of the attractor
of the difference system in

\[ \left]-1+\min_i \left\{\frac{\left(\begin{array}{c}
                                  K \\
                                  \lfloor K/2 \rfloor
                                \end{array}\right) 2^{\lfloor
K/2\rfloor}}{\alpha_i}\right\},1-\min_i
\left\{\frac{\left(\begin{array}{c}
                                  K \\
                                  \lfloor K/2 \rfloor
                                \end{array}\right) 2^{\lfloor
K/2\rfloor}}{\alpha_i}\right\}\right[. \] Then the system is not
uniformly left invertible. $\diamondsuit$
\end{Corollary}

\begin{Theorem}\label{th_uldi_uli_finite}
Fix $\mathcal U\subset\Bbb R$. Then the set of $a\in\Bbb R$ of
degree at most $K$ for which ULDI is not equivalent to ULI (in the
system (\ref{sistema})) is discrete except for possibly $2$
accumulation points given by
\[|a|=\min_{0\not=v\in\mathcal V} |v|-1.\]
Therefore, for any fixed $\delta>0$, the set of $a$ belonging to
\begin{equation}\label{eq_admissibile_a}
\left\{a\in\Bbb R: a\ algebraic\ of\ degree\ at\ most\ K,\
\left||a|-(\min_{0\not=v\in\mathcal V} |v|-1)\right|>\delta\right\}
\end{equation}
is \emph{finite}.
\end{Theorem}

\emph{Proof:} Let us show first that there exists $M>0$ such that,
if $|a|>M$, then the system (\ref{sistema}) is NOT ULI,
independently of the fact that it is ULDI or not. Indeed, the
periodic point of order $2$ (see the proof of Proposition
\ref{prop_uldi_one_step}) in the difference system given by
\[\left\{\begin{array}{c}
    ax_1+v=x_2 \\
    ax_2-v=x_1
  \end{array}\right.\]
has the solution $x_1=\frac{-v}{a+1},x_2=\frac{v}{a+1}$. As soon as
\[\min_{0\not=v\in\mathcal V}|v|=|a|+1\ \Leftrightarrow\
|a|=\min_{0\not=v\in\mathcal V}|v|-1\] this periodic point of order
$2$ lies on $]-1,1[$ and the system is not ULDI. Moreover,
\[
\lim_{|a|\rightarrow\infty} \frac{|v|}{|a|+1}= 0.
\]
This fact, together with Theorem \ref{th_notuldi_notuli} implies
that there exists $M:\ |a|>M$ implies that the system is not ULI.

Suppose now that for a particular $a\in\Bbb R$ the system
(\ref{sistema}) is ULI but not ULDI. Then it must be
\begin{equation}\label{eq_i_a}
a-\epsilon(a)<\min_{0\not=v\in\mathcal V}|v|-1<a.
\end{equation}
Moreover a fixed $\delta$ such that
\[|a-(\min_{0\not=v\in\mathcal V}|v|-1)|>\delta\]
can be supposed to exist (because an accumulation point in
$a=\min_{0\not=v\in\mathcal V}|v|-1$ is not excluded). It's now easy
to see that, once $\delta$ is fixed, there exists a
$\delta^{\prime}=\delta^{\prime}(a)$ such that
\[|a^{\prime}-a|<\delta^{\prime}\ \Rightarrow\ a^{\prime}\ does\ not\ satisfy\ (\ref{eq_i_a}).\]
This is simply because, thanks to Theorem \ref{theo_epsilon_a_alg},
the set of algebraic $a^{\prime}$ of degree at most $K$ such that
$\epsilon(a)<\delta$ is finite. Therefore $\delta^{\prime}(a)$ can
be indeed taken independently of $a$, and Theorem is thus proved.
$\diamondsuit$

\begin{Remark}
The condition given by equation (\ref{eq_admissibile_a}) is not
important from a practical point of view, since an infinitesimal
change in the quantity $\delta$, the rate of the uniform partition
$\mathcal P$, is enough to satisfy it. $\diamondsuit$
\end{Remark}

For $|a|>2$, even if Theorem \ref{th_uldi_uli_finite} doesn't work,
we have the following Theorem, that inductively construct two
initial states and two sequences of inputs that give rise to the
same output, if a particular inequality (a bit stronger than ULDI)
is satisfied.

\begin{Theorem}\label{controesempio}
Suppose that in the system (\ref{sistema}) $|a|>2$. If there exist
$u_1,u_2\in\mathcal{U},u_1\not=u_2$ such that $|u_1-u_2|<|a|$, or
equivalently if
\[\min_{0\not=v\in\mathcal V} |v|<|a|,\]
then the system is not ULI.
\end{Theorem}

\emph{Proof:} We will consider sequences of sets of type
\begin{equation}\label{eq_set_succession}\left\{\begin{array}{l}
           S_{i+1}= \left\{a(S_{i})+u(i)\right\}\bigcap\left\{a(S_{i})+u^{\prime}(i)\right\}\bigcap \mathcal{P}(i+1) \\
           S_0 = [0,1[,
         \end{array}\right.\end{equation}
where $u(i),u^{\prime}(i)\in\{u_1,u_2\}$ and
$\mathcal{P}(i+1)\in\mathcal{P}$ is chosen at each step to maximize
the measure of $S_{i+1}$.

In the sequence (\ref{eq_set_succession}) take $u(1)=u_1$,
$u^{\prime}(1)=u_2$ and $\mathcal{P}(1)$. Since $|u_1-u_2|<a$, there
exists a $\mathcal{P}(1)\in\mathcal{P}$ such that $\mu(S_1)>0.$
Then, for $i>1$ define
\[u(i)=u^{\prime}(i)=u_1.\]
Since $|a|>2$ there exists an $i_0$ such that
$\mu\left(S_{i_0}\right)=1$, therefore, applying again
$u(i_0+1)=u_1$ and $u^{\prime}(i_0+1)=u_2$
\[\mu\left\{A(S_{i_0})+Bu_1\ \cap\ A(S_{i_0})+Bu_2\right\}>0.\]
So there exists $x_0,x^{\prime}_0\in\Bbb R$ and
$(u(1),\ldots,u(i_0+1)),(u^{\prime}(1),\ldots,u^{\prime}(i_0+1)),$
with $u(1)\not=u^{\prime}(1)$ and $u(i_0+1)\not=u^{\prime}(i_0+1)$,
such that for the corresponding outputs it holds
\[(y(0),\ldots,y(i_0+1)) = (y^{\prime}(0),\ldots,y^{\prime}(i_0+1))\]

It is then enough to point out that, since we can achieve every pair
of states $x,x'\in S_{i_0}$ in the above described way, we can again
go on in the same way and find a new instant $i_1$, a pair of
initial states $x_{1,0},x'_{1,0}$, and control sequences
$(u(1),\ldots,u(i_1)),\ (u^{\prime}(1),\ldots,u'(i_1))$, with
$u(i_1)\not=u^{\prime}(i_1)$, such that for the corresponding output
it holds
\[(y(0),\ldots,y(i_1)) = (y^{\prime}(0),\ldots,y^{\prime}(i_1)).\]

Finally, we can achieve by induction an increasing finite sequence,
but arbitrarily long, of instants $i_k$, pairs of initial states
$(x_{k,0},x'_{k,0})$, and sequences of controls
$(u(1),...,u(i_k)),(u'(1),...,u'(i_k))$ with $u(i)\not =u'(i)$ if
$i=i_j+1$ for $j=1,...,k-1$ such that such that for the
corresponding output it holds
\[(y(0),\ldots,y(i_k)) = (y^{\prime}(0),\ldots,y^{\prime}(i_k)).\]
This contradicts the uniform left invertibility property.
$\diamondsuit$\\

Before giving some examples we observe that our original aim, to
show the equivalence between ULDI and ULI, has been reached, modulo
cases described in theorem \ref{th_uldi_uli_finite}. This
equivalence is actually stronger than what we showed: indeed we
didn't take into account any influence of input sequences in proofs!

\section{Examples}

\begin{Example}
Consider the system
\begin{equation}\label{eq_syst_example}
\left\{
\begin{array}{l}
  x(k+1) = ax(k)+u(k) \\
  y(k) = \lfloor x(k) \rfloor \\
  \mathcal U =
  \{-M\delta,-(M-1)\delta,\ldots,0,\ldots,(M-1)\delta,M\delta\},
\end{array}
\right.
\end{equation}
where $a,x(k),u(k),\delta >0\in\Bbb R$, $y(k)\in\Bbb Z$, $M\in\Bbb
N$. Straightforward calculations show that $\mathcal V =
\{-2M\delta,-(2M-1)\delta,\ldots,0,\ldots,(2M-1)\delta,2M\delta\}$.
For any fixed $\delta$, following the proof of Theorem
\ref{th_uldi_uli_finite}, the solutions $a$ of the equation
(\ref{eq_i_a}) should be studied:
\[a-\epsilon(a)<\delta-1<a.\]
With regard this example $a$ is supposed to be rational,
\[a=\frac{p}{q}\qquad p>q>0,\]
because it is possible to exclude the case $|a|<1$ (that can be
solved with methods described in \cite{Nev}) and because the cases
$p<0,q>0$ or $p>0,q<0$ can be obtained in a similar way. So, calling
$\tau=\delta-1$, suppose $a=\frac{p}{q}$ with $p$ of the form
\[\lfloor \tau q\rfloor +k,\qquad k\geq1.\]
So equation (\ref{eq_i_a}) becomes
\[\frac{\lfloor \tau q\rfloor +k}{q}-\tau < \frac{1}{\lfloor \tau q\rfloor +k}\quad \Leftrightarrow\]
\[\Leftrightarrow \quad \frac{k-frac(\tau q)}{q} < \frac{1}{\lfloor \tau q\rfloor +k}.\]
This last equation implies that $k=1$ (otherwise the first member
would be greater than $\frac{1}{q}$ and the second smaller than
$\frac{1}{\tau q}$). So it must be
\begin{equation}\label{eq_example}
1-frac(\tau q) < \frac{q}{\lfloor \tau q\rfloor +1}.
\end{equation}
It's obvious now that, if $\delta\not\in \Bbb Q$, since the
fractional parts
\[\left\{frac(\tau q):\ q\in\Bbb N\right\}\]
are dense in $[0,1[$, there is an infinite set of $q\in\Bbb N$ such
that (\ref{eq_example}) is satisfied, and so there exists an
infinite set of rational $a$ such that (\ref{eq_i_a}) is satisfied,
i.e. an infinite set of $a$ such that ULDI is not equivalent to ULI.
Therefore the two possible accumulation points given by
$|a|=\delta-1$ are effectively present.

Considering instead only the $a$'s belonging to the set
\[\left\{a\in\Bbb Q:\ \left|a-(\delta-1)\right|>\theta\right\},\]
the following is obtained
\[\left\{\begin{array}{l}
    \tau + \theta<a\\
    a-\epsilon(a) < \tau
    \end{array}\right.\quad\Rightarrow\quad \left\{\begin{array}{l}
                                                        a>\tau + \theta\\
                                                        \epsilon(a)>\delta
                                                        \end{array}\right.\quad\Rightarrow\quad\left\{\begin{array}{l}
                                                                                                            \frac{p}{q}>\tau + \theta\\
                                                                                                            \frac{1}{p}>\delta
                                                                                                            \end{array}\right..\]
In this case the set of $a$'s for which ULDI is not equivalent to
ULI must be found among the solutions of the latter system, and is
clearly finite: this is the set of rationals with numerator
$p<\frac{1}{\delta}$ and denominator $q<\frac{p}{\tau+\theta}$.

Suppose instead $\tau=\frac{l}{m}\in\Bbb Q$. Then (\ref{eq_example})
becomes
\[1-frac \left(\frac{l}{m} q\right) < \frac{q}{\lfloor \frac{l}{m} q\rfloor +1}.\]
In this case note that the left-hand side can assume $m$ possible
values $h_i$ for $q$ varying in $\Bbb N$, and that the right-hand
side tends to $\frac{1}{\tau+1}$ when $q$ tends to infinity. So, if
one of the $h_i$ is $<\frac{1}{\tau+1}$ there is an infinite set of
$a\in\Bbb Q$ such that ULDI and ULI are equivalent (there are the
two accumulation points), otherwise there is a finite set of
$a\in\Bbb Q$ (possibly empty) such that ULDI and ULI are equivalent
(no accumulation points). $\diamondsuit$
\end{Example}

\section{Conclusions}

In this paper we studied left invertibility of I/O quantized linear
systems of dimension $1$, and we proved that it is equivalent,
except for a finite number of cases (but there is the possibility of
having two accumulation points), to left D-invertibility, very easy
to detect (Proposition \ref{prop_uldi_one_step}). Notice that
algebraic conditions play a central role in investigation of left
invertibility of quantized systems as well in other fields when a
quantization is introduced (see for instance \cite{Bic,Chi}).

Future research will include further investigation on the
equivalence between left invertibility and left D-invertibility to
higher dimensions.

\section{Notations}

In this ``special'' section we collect all the notations used in
this paper, ordered as they appear.\\

\begin{enumerate}
  \item $frac(\cdot):\Bbb R\rightarrow\Bbb Z$: the function that
associates to each real number its fractional part: $frac(r) =
r-\lfloor r \rfloor;$

  \item $\mathcal P$: uniform partition, Definition \ref{def_uniform_partition};

  \item $f^{k_2}_{k_1}(x_0,u_1,\ldots,u_{k_2}):$ the sequence of
outputs $(y_{k_1},\ldots,y_{k_2})$ generated by the system
(\ref{sistema}) with initial condition $x_0$ and input string
$(u_1,\ldots,u_{k_2})$;

  \item $Q$: the set $\subset \Bbb R^2$ containing all pairs of states that are in the same element of the uniform partition $\mathcal
  P$, Definition \ref{def_set_q};

  \item $F_{U(k)}(X(k))$: the updating map of the $2$-dimensional
system (\ref{sistema2d});

  \item $z(k)$: state of the difference system, Definition \ref{def_difference_system};

  \item $\mathcal V$: $\mathcal U-\mathcal U$;

  \item $D^{k_2}_{k_1}(z_0,v_1,\ldots,v_{k_2}):$ the
sequence $\left(\pi_pz(k_1),\ldots,\pi_pz(k_2)\right)$ generated by
the difference system with initial condition $z_0$ and input string
$(v_1,\ldots,v_{k_2})$;

  \item $C_L$: the image mod. $1$ of the linear manifold given by
linear relations $L=\{L^j\}_1^J$ (Remark \ref{remark_kron});

  \item $R(x)$: generic polynomial, whose roots are $\rho_i$ and
degree is $q$ (Definition \ref{def_mahler_measure};

  \item $\mathfrak M$: Mahler measure (Definition \ref{def_mahler_measure};

  \item $\|R\|_\infty$: the norm given by the maximum modulus of the
coefficients of a polynomial (eq. (\ref{eq_norm_infty});

  \item $\partial$: topological boundary of a set;

  \item $\mu$: Lebesgue measure of a set;

  \item $Q^{\prime}$: the ``strip'' $(x_1,x_2)\in\Bbb R^2$ such that $|x_2-x_1|<1$, defined
in eq. (\ref{eq_set_qprime});

  \item $\varrho_i$: the line defined in equation (\ref{eq_varrho});

  \item $k_i$: the distance between the point $\left(\begin{array}{c}
                                                        c_i \\
                                                        c^{\prime}_i \\
                                                    \end{array}\right)$,
defined in (\ref{eq_varrho}), and the union of positive coordinate
axes, along the line $\varrho_i$ (refer to the Fig. $1$);

  \item $\Omega$: the set defined in (\ref{eq_set_omega});

  \item $P_a(x)=\alpha_0+\ldots+\alpha_Kx^K$: minimal polynomial of
  $a$, with coefficients $\alpha_i$ and degree $K$;
  \item $P_K^J,\Psi_K^J$: the linear manifold and the matrix defined
  by equations (\ref{eq_pkj}) and (\ref{eq_psikj});
\end{enumerate}


\begin{thebibliography}{99}
\bibitem{Bak}~Baker A., \emph{Trascendental number theory}, Cambridge University
Press, (1993).

\bibitem{Bar}~Barnsley M., \emph{Fractals everywhere}, Academic
Press inc, (1993).

\bibitem{Bic}~Bicchi A., Marigo A., Piccoli B., \emph{On the reachability of quantized control
sytems}. IEEE Transactions on Automatic Control, 47(4), pages:
546--563, (1992).

\bibitem{Bob}~Bobylev N.A., Emel'yanov S.V., Korovin S.K., \emph{Attractor of discrete controlled systems in metric
spaces}, Computational Mathematics and Modeling, 11(4), pages:
321--326, (2000).

\bibitem{Brok}~Brockett R.W., and Mesarovic M.D., \emph{The reproducibility of multivariable
control systems}, Journal of Mathenatical Analalysis and
Applications, 11, pages: 548--563, (1965).

\bibitem{Car}~Carli R., Fagnani F., Speranzon A., and Zampieri S.,
\emph{Communication constraints in the state agreement problem},
Automatica, (to appear).

\bibitem{Cas2}~Cassels J.W.S., \emph{An introduction to the geometry of numbers}. Springer, (1997).

\bibitem{Chi}~Chitour Y., Piccoli B., \emph{Controllability for
discrete systems with a finite control set}, Mathematics of Control
Signal and Systems 14, pages: 173--193, (2001).

\bibitem{Del}~Delchamps D.F., \emph{Stabilizing a linear system with quantized state
feedback}, IEEE Transactions on Automatic Control, 35(8), pages:
916--924, (1990).

\bibitem{Nev}~Dubbini N., Piccoli B., Bicchi A., \emph{Left invertibility of discrete systems with finite inputs and
quantized output}, Proceedings of 47--th IEEE Conference on Decision
and Control, pages: 4687--4692, (2008).

\bibitem{Dub}~Dubbini N., Piccoli B., Bicchi A., \emph{Left invertibility of
discrete systems with finite inputs and quantized output},
International Journal Of Control, note: accepted, (2009).

\bibitem{Dub1}~Dubbini N., Piccoli B., Bicchi A., \emph{Left invertibility of discrete--time I/O quantized
linear systems}, Mathematics of control signals and systems, note:
submitted, (2009).

\bibitem{Dub2}~Dubbini N., Monge M., \emph{An equivalent of Kronecker's theorem
for powers of an algebraic numbers}, arXiv:0910.5182v1 [math.NT],
(2009).

\bibitem{Edel}~Edelmayer A., Bokor J., Szab\'{o} Z., Szigeti F., \emph{Input reconstruction by means of system
inversion: a geometric approach to fault detection and isolation in
nonlinear systems}, International journal of applied mathematics and
computer science, 14(2), pages: 189--199, (2004).

\bibitem{FG}~Falconer K., \emph{Fractal geometry, mathematical foundations and
applications}, John Wiley and Sons, (2003).

\bibitem{Gw}~Gwozdz-Lukawska G., Jachymski J., \emph{The Hutchinson-Barnsely
theory for infinite iterated function systems}, Bulletin of
Australian Mathematical Society, 72, pages: 441--454, (2005).

\bibitem{Har}~Hardy G.H., Wright E.M.,  \emph{An introduction to the theory of
numbers}, Oxford Science Publications, (1979).

\bibitem{Inou}~Inoue E., Ushio T., \emph{Chaos communication using unknown input
observer}, Electronic and Comunication in Japan, Part 3, 84, (2001).

\bibitem{Mas}~Massey J.L., Sain M.K., \emph{Invertibility of linear time--invariant dynamical
systems}, IEEE Transactions on Automatic Control, AC-14(2), pages:
141--149, (1969).

\bibitem{Mas2}~Massey J.L., Sain M.K., \emph{Inverses of linear sequential
circuits}, IEEE Transactions on Computers, C-17, pages: 330--337,
(1968).

\bibitem{Mors}~Morse A.S., Wonham W.M., \emph{Status of noninteracting control}, IEEE
Transactions on Automatic Control, 16(6), 568--581, (1971).

\bibitem{Pica}~Picasso B., Bicchi A., \emph{On the stabilization of linear systems under assigned I/O
quantization}, IEEE Transactions on Automatic Control, 52(10),
pages: 1994--2000, (2007).

\bibitem{Resp}~Respondek W., \emph{Right and Left Invertibility of Nonlinear Control Systems},
Nonlinear Controllability and Optimal Control, New York, pages:
133--176, (1990).

\bibitem{Silv}~Silverman L.M., \emph{Inversion of multivariable linear systems}, IEEE
Transactions on Automatic Control, 14(3), pages: 270--276, (1969).

\bibitem{Son}~Sontag E.D., \emph{Mathematical control
theory: deterministic finite dimensional systems}, Springer, New
York, (1998).

\bibitem{Szn}~Szanier M., Sideris A., \emph{Feedback control of quantized constrained systems with
applications to neuromorphic controller design}, IEEE Transactions
on Automatic Control, 39(7), pages: 1497--1502, (1994).

\bibitem{Tat}~Tatikonda S.C., Mitter S., \emph{Control under communication constraints,} IEEE Transactions
on Automatic Control, 49(7), pages: 1056--1068, (2004).

\bibitem{Vu}~Vu L., Liberzon D., \emph{Invertibility of switched linear
systems}, Proceedings of the 45th IEEE Conference on Decision and
Control, pages: 4081--4086, (2006).
\end{thebibliography}
\end{document}